\documentclass[11pt,draft]{amsart}

\newtheorem{theorem}{Theorem}[section]

\theoremstyle{definition}

\numberwithin{equation}{section}

\RequirePackage{amsmath}
\RequirePackage{amssymb}
\RequirePackage{amsthm}
 \RequirePackage{epsfig}

\begin{document}

\date{}

\title[Multi-fold sums from a set with few products]
{Multi-fold sums from a set with few products}

\author{Liangpan Li}

\address{Department of Mathematical Sciences,
Loughborough University, LE11 3TU, UK }
 \email{liliangpan@gmail.com}

 \renewcommand{\thefootnote}{}




\date{}

\begin{abstract}
In this paper we show that for any $k\geq2$, there exist two
universal constants $C_k,D_k>0$, such that for any finite subset $A$
of positive real numbers with $|AA|\leq M|A|$, $|kA|\geq
\frac{C_k}{M^{D_k}}\cdot|A|^{\log_42k}.$
\end{abstract}

\maketitle

\section{Introduction}

We begin with some notation: Given a finite subset $A$ of some
commutative ring, we let $A\star A$ denote the set $\{a\star
b:a,b\in A\}$, where $\star$ is a binary operation on $A$. When
three or more summands or multiplicands are used, we let $kA$ denote
the $k$-fold sum-set $A+A+\cdots+A$, and let $A^{(k)}$ denote the
$k$-fold product-set $AA\cdots A$.

Erd\H{o}s and Szemer\'{e}di (\cite{Erdos}) once conjectured that for
any $\alpha<2$, there exists a universal constant $C_{\alpha}>0$,
such that for finite subset $A$ of real numbers,
\[\max\{|A+A|,\  |AA|\}\geq C_{\alpha}|A|^{\alpha}.\]
Non-trivial lower bounds for $\alpha$ were achieved by many authors
such as by Erd\H{o}s and Szemer\'{e}di (\cite{Erdos},
qualitatively), Nathanson (\cite{Nathanson}, $32/31$), Ford
(\cite{Ford}, $16/15$), Chen (\cite{Chen}, $6/5$), Elekes
(\cite{Elekes}, $5/4$), and Solymosi (\cite{S2}, $14/11-o(1)$;
\cite{Solymosi}, $4/3-o(1)$).

Another type of question than one can attack regarding sums and
products is to either assume that the sum-set $A+A$ is very small,
and then to show that the product-set $AA$ is very large, or to
suppose that $AA$ is very small, and then to show that $A+A$ is very
large. The best two results toward this question are respectively
due to Elekes and Ruzsa (\cite{ElekesRuzsa}), who fully confirmed
the first part of the question, and by Chang (\cite{Chang}), who
solved the second part of the question in the setting of
\textsf{integers}.

Similarly, one can consider multi-fold sums and products, but very
few results are known especially in the setting of reals. Let $B$ be
a finite subset of \textsf{integers}, then Chang (\cite{Chang})
showed that if $|BB|\leq|B|^{1+\epsilon}$, then the multi-fold
sum-set $|kB|\gg_{\epsilon,k}|B|^{n-\delta}$, where
$\delta\rightarrow0$ as $\epsilon\rightarrow0$; and Bourgain and
Chang (\cite{BourgainChang}) proved that for any $b\geq1$, there
exists $k\in\mathbb{N}$ independent of $B$ such that
$|kB|\cdot|B^{(k)}|\geq|B|^b.$ At the moment how to extent these
results to the real numbers is not known yet. Recently, Croot and
Hart established in \cite{CrootHart} the following interested
result:
\begin{theorem}\label{CrootHart}
For all $k\geq2$ and $\epsilon\in(0,\epsilon_0(k))$ we have that the
following property holds for all $n>n_0(k,\epsilon)$: If $A$ is a
set of $n$ real numbers and $|AA|\leq n^{1+\epsilon}$, then
\[|kA|\geq n^{\log_4k-f_k(\epsilon)},\]
where $f_k(\epsilon)\rightarrow0$ as $\epsilon\rightarrow0$.
\end{theorem}

 Croot and Hart also remarked that they have several different approaches to proving a theorem of the quality of  Theorem \ref{CrootHart}.

 The purpose of the present paper is to give the following slight improvement of the above Croot-Hart theorem in a rather elementary
way.  Our idea comes from Solymosi's wonderful proof
(\cite{Solymosi}) of the best currently known sum-product estimates
of real numbers mentioned earlier. Solymosi's idea has appeared
elsewhere in \cite{Hart} and \cite{LiShen}.

\begin{theorem}\label{MainTheorem}
For any $k\geq2$, there exist three positive universal constants
$C_k$, $D_k$, $\Psi_k$, such that for any finite subset $A$ of
positive real numbers with $|AA|\leq M|A|$,
\[|kA|\geq \frac{C_k}{M^{D_k}}\cdot|A|^{\Psi_k}.\]
With $\Psi_1\triangleq1$, the constants $\{\Psi_k\}_{k\geq2}$ can be
generated in any of the following way:
\[\Psi_k=\frac{1+\Psi_{k_1}+\Psi_{k_2}}{2}\ \
(k_1+k_2=k).\] Particularly, we can take $\Psi_k=\log_42k$.
\end{theorem}


There are some other interested estimates on sum-sets and
product-sets in the reals. For example, see
\cite{ElekesNathansonRuzsa}, \cite{IosevichRocheRudnev},
\cite{Schoen} and \cite{Shen}.



\section{Proof of the main theorem}

We will prove Theorem \ref{MainTheorem} for all $k\in \mathbb{N}$ by
induction. Obviously, one can choose $D_1=0$, $C_1=\Psi_1=1$. Next
for any $k\geq2$, we assume the existences of positive universal
constants $C_i$, $D_i$ and $\Psi_i$ for all $i\in[2,k)$. Our purpose
is to find $C_k$, $D_k$ and $\Psi_k$ satisfying the required
property. Let $k_1,k_2$ be any two natural numbers such that
$k_1+k_2=k$.

By the Ruzsa triangle inequality, $|A/A|\leq M^2|A|$. For any $s\in
A/A$, let $A_s\triangleq \{(x,y)\in A\times A: y=sx\}$. Let
$D=\{s:|A_s|\geq\frac{|A|}{2M^2}\}$, and let $s_1<s_2<\cdots<s_m$
denote the elements of $D$, labeled in increasing order. Obviously,
\[\sum_{s\in D}|A_s|\geq\frac{|A|^2}{2},\]
which implies $m\geq\frac{|A|}{2}$. Let $A_{m+1}$ be the projection
of $A_m$ onto the vertical line $x=\min A$, and let
$\Pi:\mathbb{R}^2\rightarrow \mathbb{R}$ be the projection map from
$\mathbb{R}^2$ onto the vertical axis.  It is geometrically evident
that $\{k_1A_j+k_2A_{j+1}\}_{j=1}^m$ are mutually disjoint.
 Thus
\[|(kA)\times(kA)|\geq\sum_{j=1}^m|k_1A_j+k_2A_{j+1}|=\sum_{j=1}^m|k_1A_j|\cdot|k_2A_{j+1}|=\sum_{j=1}^m|k_1\Pi(A_j)|\cdot|k_2\Pi(A_{j+1})|.\]
Note
\[|\Pi(A_j)\Pi(A_j)|\leq|AA|\leq M|A|\leq 2M^3|\Pi(A_j)|.\]
Applying induction to all of the $\Pi(A_j)$'s,
\[|kA|^2\geq \frac{|A|}{2}\cdot \frac{C_{k_1}}{(2M^3)^{D_{k_1}}}(\frac{|A|}{2M^2})^{\Psi_{k_1}}\cdot \frac{C_{k_2}}{(2M^3)^{D_{k_2}}}(\frac{|A|}{2M^2})^{\Psi_{k_2}},\]
which yields
\[|kA|\geq\Big(\frac{C_{k_1}\cdot C_{k_2}}{2\cdot (2M^3)^{D_{k_1}+D_{k_2}}\cdot(2M^2)^{\Psi_{k_1}+\Psi_{k_2}}}\Big)^{1/2}\cdot|A|^{\frac{1+\Psi_{k_1}+\Psi_{k_2}}{2}}.\]
Thus one can let $\Psi_k\triangleq\frac{1+\Psi_{k_1}+\Psi_{k_2}}{2}$
and define $C_k, D_k$ in a similar way.

Finally, let $z\triangleq\lfloor\log_2k\rfloor$. Then
\[\Psi_{2^z}\geq\frac{1}{2}+\Psi_{2^{z-1}}\geq\cdots\geq\frac{z}{2}+\Psi_1=\frac{z+2}{2}\geq\log_42k.\]
Consequently,
\[|kA|\geq|2^zA|\geq\frac{C_{2^z}}{M^{D_{2^z}}}\cdot|A|^{\Psi_{2^z}}\geq\frac{C_{2^z}}{M^{D_{2^z}}}\cdot|A|^{\log_42k}.\] This
concludes the whole proof.

\textbf{Acknowledgements.} The author was supported by the NSF of
China (11001174).


\begin{thebibliography}{99}


\bibitem{BourgainChang}
J. Bourgain, M.-C. Chang, \emph{On the size of $k$-fold sum and
product sets of integers}, J. Amer. Math. Soc. 17 (2003) 473--497.

\bibitem{Chang}
M.-C. Chang, \emph{The Erd\H{o}s-Szemer\'{e}di problem on sum set
and product set}, Ann. Math. 157 (2003) 939--957.


\bibitem{Chen}
Y.~G. Chen, \emph{On sums and products of integers}, Proc. Amer.
Math. Soc. 127 (1999) 1927--1933.



\bibitem{CrootHart}
E. Croot, D. Hart, \emph{$h$-fold sums from a set with few
products}, SIAM J. Discrete Math. 24 (2010) 505--519.

\bibitem{Elekes}
Gy. Elekes, \emph{On the number of sums and products}, Acta Arith.
81 (1997) 365--367.

\bibitem{ElekesNathansonRuzsa}
Gy. Elekes, M.~B. Nathanson, I.~Z. Ruzsa, \emph{Convexity and
sumsets}, J. Number Theory 83 (1999) 194--201.

\bibitem{ElekesRuzsa}
Gy. Elekes, I.~Z. Ruzsa, \emph{Few sums, many products}. Studia Sci.
Math. Hungar. 40 (2003) 301--308.



\bibitem{Erdos}
P. Erd\"{o}s and E. Szemer\'{e}di, \emph{On sums and products of
integers}. In: Studies in Pure Mathematics (Birkhauser, Basel, 1983)
213--218.

\bibitem{Ford}
K. Ford, \emph{Sums and products from a finite set of real numbers},
Ramanujan J. 2 (1998) 59--66.

\bibitem{GaraevShen}
M.~Z. Garaev, C.-Y. Shen, \emph{On the size of the set $A(A+1)$},
Math. Z. 265 (2010) 125--132.

\bibitem{Hart}
D. Hart, A. Niziolek, \emph{Some results on the size of sum and
product sets of finite sets of real numbers}, Involve 2 (2009)
603--609.


\bibitem{IosevichRocheRudnev}

A. Iosevich, O. Roche-Newton, M. Rudnev, \emph{On an application of
Guth-Katz theorem}, arXiv:1103.1354, accepted by Math. Research
Letters, 2011.

\bibitem{LiShen}
L. Li, J. Shen, \emph{A sum-division estimate of reals}, Proc. Amer.
Math. Soc. 138 (2010) 101--104.


\bibitem{Nathanson}
M.~B. Nathanson, \emph{On sums and products of integers}, Proc.
Amer. Math. Soc. 125 (1997) 9--16.

\bibitem{Schoen}
T. Schoen, I.~D. Shkredov, \emph{On sumsets of convex sets},
arXiv:1105.3542, 2011.

\bibitem{Shen}
C.-Y. Shen,  \emph{Algebraic methods in sum-product phenomena},
arXiv:0911.2627, to appear in Israel J. Math., 2009.


\bibitem{S2}
J. Solymosi, \emph{On the number of sums and products}, Bull. London
Math. Soc. 37 (2005) 491--494.

\bibitem{Solymosi} J. Solymosi, \emph{Bounding multiplicative energy by the
sumset}, Adv. Math. 222 (2009) 402--408.




















\end{thebibliography}
\end{document}